\documentclass{article}
\usepackage{float}

\usepackage[final]{neurips_2025_ml4ps}

\usepackage[utf8]{inputenc}
\usepackage[T1]{fontenc}
\usepackage{hyperref}
\usepackage{url}
\usepackage{booktabs}
\usepackage{amsfonts}
\usepackage{nicefrac}
\usepackage{microtype}
\usepackage{xcolor}
\usepackage{amsmath}

\usepackage{graphicx}

\title{One-Shot Transfer Learning for Nonlinear PDEs with Perturbative PINNs}

\author{%
  \textbf{Samuel Auroy} \\
  \normalfont ENSTA Paris \\
  \normalfont\texttt{samuel.auroy@ensta.fr}
  \and
  \textbf{Pavlos Protopapas} \\
  \normalfont Harvard University \\
  \normalfont\texttt{pprotopapas@g.harvard.edu}
}

\begin{document}

\maketitle

\begin{abstract}
We propose a framework for solving nonlinear partial differential equations (PDEs) by combining perturbation theory with one-shot transfer learning in Physics-Informed Neural Networks (PINNs). Nonlinear PDEs with polynomial terms are decomposed into a sequence of linear subproblems, which are efficiently solved using a Multi-Head PINN. Once the latent representation of the linear operator is learned, solutions to new PDE instances with varying perturbations, forcing terms, or boundary/initial conditions can be obtained in closed form without retraining. 

We validate the method on KPP–Fisher and wave equations, achieving errors on the order of $10^{-3}$ while adapting to new problem instances in under $0.2$ seconds; comparable accuracy to classical solvers but with faster transfer. Sensitivity analyses show predictable error growth with $\epsilon$ and polynomial degree, clarifying the method’s effective regime. 

Our contributions are: (i) extending one-shot transfer learning from nonlinear ODEs to PDEs; (ii) deriving a closed-form solution for adapting to new PDE instances; and (iii) demonstrating accuracy and efficiency on canonical nonlinear PDEs. We conclude by outlining extensions to derivative-dependent nonlinearities and higher-dimensional PDEs.
\end{abstract}

\section{Introduction}
Partial differential equations (PDEs) are central to modeling in biology, physics, and engineering. Unlike ordinary differential equations (ODEs), which govern single-variable dynamics, PDEs describe coupled spatial-temporal processes and are harder to solve due to stability, dimensionality, and complex boundaries. Classical ODE methods like Runge--Kutta \citep{BUTCHER1996247} extend only partially, leaving finite element and related methods \citep{ciarlet:hal-04039611} as workhorses, but at high computational cost.

Recent machine learning approaches, especially Physics-Informed Neural Networks (PINNs), provide an appealing alternative \citep{Lagaris_1998, wong2025evolutionaryoptimizationphysicsinformedneural,RAISSI2019686,SIRIGNANO20181339,krishnapriyan2021characterizingpossiblefailuremodes}. By embedding physical laws into neural networks, PINNs bypass discretization and yield differentiable solutions, but each new instance typically requires retraining. Multi-instance strategies such as bundle training \citep{flamant2020solvingdifferentialequationsusing} or Multi-Head PINNs \citep{zou2023lhydramultiheadphysicsinformedneural} help, yet generalization across PDE families remains limited.

One-shot transfer learning shows strong results for linear equations \citep{desai2022oneshottransferlearningphysicsinformed}, but nonlinear cases still need iterative optimization. Building on \citep{lei2023oneshottransferlearningnonlinear}, which addressed nonlinear ODEs with perturbative terms, we extend this approach to PDEs with polynomial perturbations. We validate it on KPP--Fisher equations and argue for broader applicability across nonlinear PDEs. Extending to higher-dimensional PDEs (for example, Navier--Stokes) is future work; here we validate on canonical PDEs with well-understood dynamics to isolate the method's behavior.
\newpage
\paragraph{Contributions.} 
This work makes the following contributions:  
\begin{itemize}
    \item We extend one-shot transfer learning, previously demonstrated only for nonlinear ODEs, to nonlinear PDEs with polynomial perturbations.  
    \item We derive a perturbative framework that reduces nonlinear PDEs into a sequence of linear subproblems, enabling efficient reuse of a shared latent representation learned by a Multi-Head PINN.  
    \item We demonstrate the method on canonical PDEs (KPP--Fisher and wave equations), showing accuracy comparable to classical solvers while adapting more efficiently to new forcing terms, boundary conditions, and perturbations.  
    \item We analyze sensitivity to perturbation size and polynomial degree, and discuss extensions to related operators and non-perturbative terms.  
\end{itemize}

\section{Methodology}
We study a class of scalar nonlinear PDEs of order $n$,
\begin{equation}
\mathcal{D} u + \epsilon P(u) = f(x,t),
\label{eq:fullproblem}
\end{equation}
\begin{equation}
\text{where}\quad \mathcal{D} u = \sum_{|\alpha|=0}^{n} a_{\alpha}(x,t)\,\partial^{\alpha}u \ \text{for } \alpha \in \mathbb{N}^2,\quad
P(u) = \sum_{l=0}^{m} P_l u^l,\ \text{and}\ \forall l\in [0,m],\ P_l \leq 1.
\label{eq:def}
\end{equation}
Here, $\epsilon < 1$ ensures that the polynomial term $\epsilon P(u)$ acts as a perturbation to the linear PDE defined by $\mathcal{D}$. Because all coefficients $P_l$ are bounded by 1, each nonlinear contribution $\epsilon P_l u^l$ remains small relative to the linear term. Appropriate boundary and initial conditions (Dirichlet or Neumann) are imposed so that the problem is well-posed.

\subsection{Expansion of the solution}
\label{expansion}
Intuitively, we separate the PDE into a dominant linear operator $\mathcal{D}$ and a smaller nonlinear correction $\epsilon P(u)$. This allows us to write the solution as a series of corrections of increasing order in $\epsilon$. Each correction is governed by a linear PDE, so the nonlinear problem is reduced to a sequence of linear problems.

Following \citep{lei2023oneshottransferlearningnonlinear} with the technique from \citep{Pertubation_Nayfeh}, assume $u(x,t)= \sum_{i=0}^{\infty} \epsilon^{i} u_i(x,t)$ and truncate,
\begin{equation}
u(x,t) \approx \sum_{i=0}^{p} \epsilon^{i} u_i(x,t).
\end{equation}
For this perturbative expansion to be meaningful, we assume $\epsilon < 1$, so higher-order terms diminish.

Substituting into \eqref{eq:fullproblem} gives
\begin{equation}
\sum_{i=0}^{p} \epsilon^i \mathcal{D} u_i + \epsilon \sum_{l=0}^{m} P_l
\left( \sum_{\substack{k_0 + k_1 + \ldots + k_p = l}}
\frac{l!}{k_1! k_2! \ldots k_p!}
\epsilon^{\sum_{i=0}^{p} i k_i}
\prod_{i=0}^{p} u_i^{k_i} \right) = f.
\end{equation}
The left-hand side is a polynomial in $\epsilon$. Collecting equal powers yields $p+1$ linear PDEs.

At order $\epsilon^0$,
\begin{equation}
\mathcal{D} u_0 = f.
\end{equation}
For a general order $j$ with $1 \leq j \leq p$,
\begin{equation}
\mathcal{D} u_j = - \sum_{l=0}^{m} P_l
\left( \sum_{\substack{k_0 + k_1 + \ldots + k_p = l \\
\sum_{i=0}^{p} i k_i = j - 1}}
\frac{l!}{k_1! k_2! \ldots k_p!}
\prod_{i=0}^{p} u_i^{k_i} \right) := f_j.
\end{equation}
Thus we obtain linear PDEs of the form
\begin{equation}
\mathcal{D} u_j(x,t) = f_j(x,t),
\end{equation}
where each $f_j$ depends only on $u_i$ with $i < j$. We solve them iteratively to construct
\begin{equation}
u(x,t) \approx \sum_{i=0}^{p} \epsilon^{i} u_i(x,t).
\end{equation}
To satisfy ICs and BCs, a common choice is to impose the full problem's ICs/BCs on $\mathcal{D} u_0 = f$ and homogeneous ICs/BCs on higher orders.

\subsection{Multi-Head Physics-Informed Neural Network (MH-PINN)}
We solve the linear PDEs $\mathcal{D} u_j = f_j$ using a Multi-Head PINN. The network is fully connected with $K$ heads. The $k$th output is
\begin{equation}
u_k(x,t) = H(x,t) W_k,
\end{equation}
where $H(x,t)$ is the last shared layer's activation and $W_k$ are head-specific weights.

We compute a physics-informed loss for each head and average:
\begin{align}
\mathcal{L} &= \tfrac{1}{K}\sum_{k=0}^{K}\mathcal{L}_{k}, \\
\mathcal{L}_k &= w_{pde}\mathcal{L}_{k,\text{pde}} 
+ w_{IC}\mathcal{L}_{k,\text{IC}} 
+ w_{BC}\mathcal{L}_{k,\text{BCs}} \\
\nonumber&= w_{pde}\left( \mathcal{D} u_k(x,t) - f_k(x,t) \right)^2 
+ w_{IC}\left( u_k(x,0) - g_k(x) \right)^2
+ w_{BC}\sum_{\mu \in \{L,R\}} \left( u_k(\mu,t) - B_{\mu,k}(t) \right)^2.
\end{align}
Here $w_{pde}, w_{IC}, w_{BC} > 0$ balance PDE residual, initial conditions, and boundary conditions. Neumann conditions are handled by replacing $u_k$ with $\partial^\alpha u_k$ in IC/BC terms. The total loss is computed over collocation points and minimized by backpropagation. The key objective is to learn a reusable latent space $H$ tied to $\mathcal{D}$.

\subsection{One-shot Transfer Learning}
Let $(f^*, g^*, B^*)$ denote the forcing, initial data, and boundary data for a new linear problem. With the body frozen, evaluate $H$ at sampled points $(\hat{x}, \hat{t})$ and write
\begin{equation}
u^*(\hat{x},\hat{t}) = H(\hat{x},\hat{t}) W^*.
\end{equation}
Define the loss
\begin{equation}
\mathcal{L}= w_{pde}\left( \mathcal{D} H W^* - f^*\right)^2 
+ w_{IC}\left( H_0 W^* - g^* \right)^2
+ w_{BC}\sum_{\mu \in \{L,R\}} \left( H_\mu W^* - B^*_{\mu}\right)^2,
\end{equation}
where $H_0$ and $H_\mu$ denote $H$ at initial and boundary points. Each term is convex in $W^*$, so $\partial \mathcal{L}/\partial W^* = 0$ yields
\begin{equation}
    W^* = M^{-1} \left( 
    w_{pde}\mathcal{D}_H^T f^* + 
    w_{IC} H_0^T g^* +  
    w_{BC}\sum_{\mu \in \{L,R\}} H_\mu^T B_\mu^* 
    \right),
\end{equation}
with
\begin{equation}
M = w_{pde}\mathcal{D}_H^T \mathcal{D}_H 
+ w_{IC} H_0^T H_0 
+ w_{BC}\sum_{\mu \in \{L,R\}} H_\mu^T H_\mu.
\end{equation}
Importantly, $M$ depends only on $\mathcal{D}$ and the sampling strategy, not on $(f^*, g^*, B^*)$, so $M^{-1}$ can be reused across new tasks.

\section{Results}
We applied our methodology to a family of KPP Fisher equations (details on MH-PINN training are provided in Section~\ref{subsec: training}):
\begin{equation}
    \frac{\partial u}{\partial t} - D \frac{\partial^2 u}{\partial x^2} - \epsilon u^{n_1} (1 - u^{n_2}) = f(x,t).
\end{equation}
This is the form in \eqref{eq:fullproblem} with $\mathcal{D} = \frac{\partial}{\partial t} - D \frac{\partial^2 }{\partial x^2}$ and $P(u)= -u^{n_1} +u^{n_1+n_2}$. The equation first appeared in population dynamics \citep{NADIN2012633eqorigin,ALFARO20171309eqmodel} and has been applied in chemistry \citep{brunet:tel-01417420}. When $f=0$, it describes evolution between equilibria 0 and 1. For $n_1 = n_2 = 1$ and the conditions from Table~\ref{table: equations}, we compute a solution under conditions similar to those in \citep{kppTrefethen} and observe the propagation toward the equilibrium value 1.

\begin{figure}[H]
    \centering
    \begin{minipage}{0.7\textwidth}
        \centering
        \includegraphics[width=\linewidth]{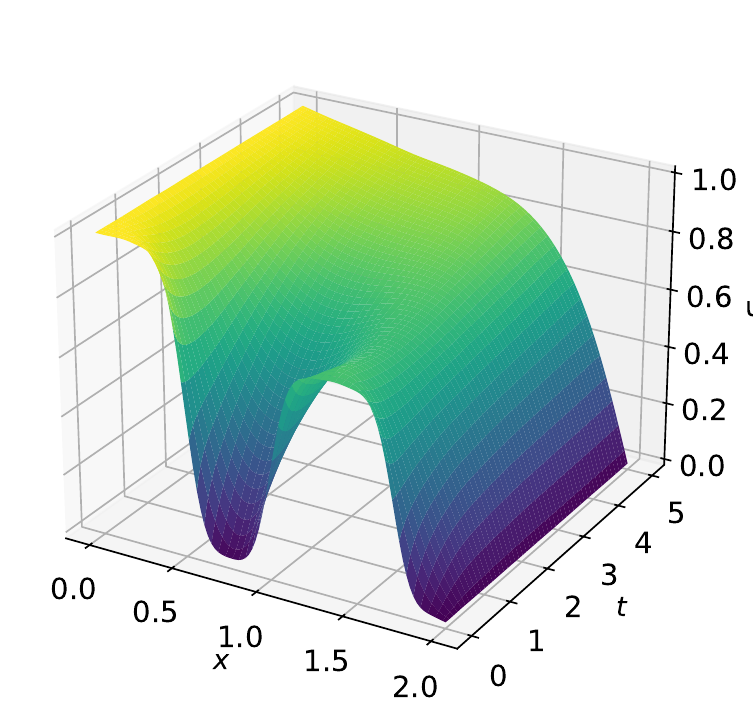}
        \caption{Solution to the nonlinear equation computed by PINNs.}
        \label{fig:1st example}
    \end{minipage}
\end{figure}
 We then compare the error obtained with our approach to that of an approximate solution computed using a classical numerical solver, as a function of the perturbation order $p$.
\begin{figure}[H]
    \centering
    \begin{minipage}{0.7\textwidth}
        \centering
        \includegraphics[width=\linewidth]{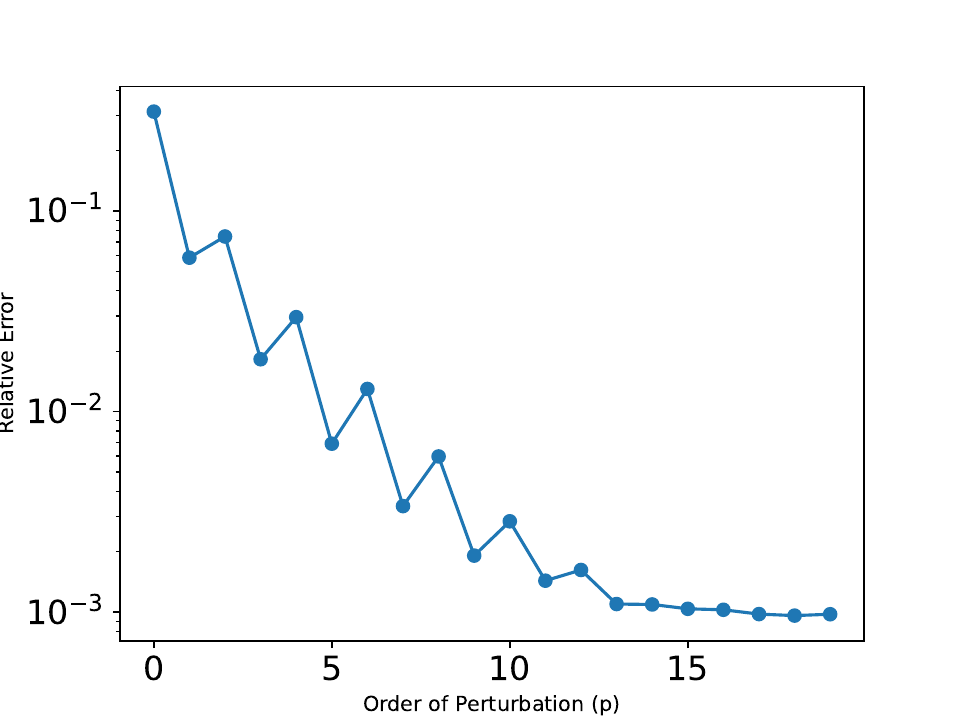}
        \caption{Relative error versus perturbation order $p$.}
        \label{fig:error for pertubation}
    \end{minipage}
\end{figure}

Our method solves this equation in $1.489 \cdot10^{-1}$ seconds for $p = 10$, compared to $1.62 \cdot10^{-1}$ seconds using a method-of-lines solver with SciPy's \texttt{solve\_ivp}.

The solution can be computed for various $n_1$ and $n_2$. Increasing $n_2$ accelerates propagation toward equilibrium 1, while larger $n_1$ slows it. We observe similar precision and runtime for the fourth set of conditions in Table~\ref{table: equations}. The relative errors are $1.1 \cdot 10^{-3}$, $1.2 \cdot 10^{-3}$, and $1.9 \cdot 10^{-2}$ respectively; plots are in Appendix \ref{Graph solutions}. The slightly higher error in the third case stems from its larger $\epsilon$ (see Section \ref{Limitations}).

To test robustness, we trained another Multi-Head PINN for the wave operator $\mathcal{D} = \frac{\partial^2}{\partial t^2} - c^2 \frac{\partial^2 }{\partial x^2}$ and obtained similar results (Appendix \ref{Results wave equation}).

\section{Limitations and Extensions}
\label{Limitations}

The method's accuracy is sensitive to the choice of $\epsilon$. As shown in Appendix~\ref{Graph solutions}, the error exhibits a sharp increase once $\epsilon$ exceeds a threshold that depends on solution amplitudes and forcing. Specifically, solutions with larger amplitudes require smaller $\epsilon$ values to maintain accuracy. A favorable aspect of the method is that the problematic range of $\epsilon$ values is easily identifiable in experiments, as solutions diverge when $\epsilon$ becomes too large. For higher-degree polynomials (degree $>10$), even smaller $\epsilon$ values are required due to the rapid growth of perturbation terms.

This work lays the foundation for numerous extensions:
\begin{itemize}
    \item The approach can transfer across related differential operators, for instance when the coefficients $a_\alpha(t)$ in Eq.~\ref{eq:def} vary slightly, as demonstrated in \citep{seiler2025stifftransferlearningphysicsinformed,lei2023oneshottransferlearningnonlinear}. 

    \item The derivation naturally extends to polynomial terms involving $u$ and its derivatives. In this case, the mathematical extension is straightforward: one defines the multi-variable polynomial and collects terms for each power of epsilon. The functions $f_j$ now depend on previously computed solutions and their derivatives, which represents the additional requirement for the experimental setup. Testing this approach on perturbative Burgers equations constitutes a natural next step.

\end{itemize}

 A key open challenge remains extending one-shot transfer learning to non-perturbative terms, which will likely require the development of new mathematical tools.

\section{Conclusion}
We developed an efficient method for solving nonlinear PDEs, extending the approach of \citep{lei2023oneshottransferlearningnonlinear} for ODEs. By combining perturbation theory with one-shot transfer learning, our method attains accuracy comparable to an optimized classical solver while adapting faster to new instances, positioning it between single-task PINNs and full-scale operator learners (e.g. Neural Operators in \citep{Neural_operator,li2021fourierneuraloperatorparametric,Lu_2021}), our model generalizes to diverse cases with a leaner backbone and minimal fine-tuning. As discussed in Section \ref{Limitations}, several extensions are promising, including derivative-dependent nonlinearities and mild operator changes, which could broaden applicability to larger-scale physical systems.
\newpage
\section{Appendix}

The code used to generate the results of this study is publicly available on this GitHub: https://github.com/Sam147258369/NeurIPS-2025-Submission-Code

\subsection{Equations used}
We define a function used as an initial condition:
\[
h(x) = \frac{1}{1+e^{20(x-0.5)}} + 0.7\left(\frac{1}{1+e^{-20(x-1)}} + \frac{1}{1+e^{-20(x-1.5)}}\right),
\]
a smooth approximation of a discontinuous function equal to 1 for $x<0.5$, equal to 0.7 for $1<x<1.5$, and 0 elsewhere. This is the same type of condition as in \citep{kppTrefethen}, but PINNs struggle to learn discontinuous functions.

\label{table: equations}
\begin{tabular}{|c|c|c|c|c|c|}
\hline
Figure & $\epsilon$ & IC & BC Left & BC Right  & Forcing function \\
\hline
\ref{fig:1st example} & 0.5 & $\psi(0, x) = h(x)$ & $\psi(t, 0) = 1$ & $\psi(t, 2) = 0$ & $f(t, x) = 0$ \\
\hline
\ref{fig:big n1}  & 0.1 & $\psi(0, x) = h(x)$ & $\psi(t, 0) = 1$ & $\psi(t, 2) = 0$ & $f(t, x) = 0$ \\
\hline
\ref{fig:big n2}  & 0.1 & $\psi(0, x) = h(x)$ & $\psi(t, 0) = 1$ & $\psi(t, 2) = 0$ & $f(t, x) = 0$ \\
\hline
\ref{fig:various}  & 0.3 & $\psi(0, x) = 1+0.5x\sin(2\pi x)$ & $\psi(t, 0) = e^{-t}$ & $\psi(t, 2) = 1 + \sin(t)$ & $f(t, x) = \sin(2t)\cos(3x)$ \\
\hline
\end{tabular}

\subsection{Graphs of various solutions}
\label{Graph solutions}

\begin{figure}[ht]
    \centering
    \begin{minipage}{0.32\textwidth}
        \centering
        \includegraphics[width=\linewidth]{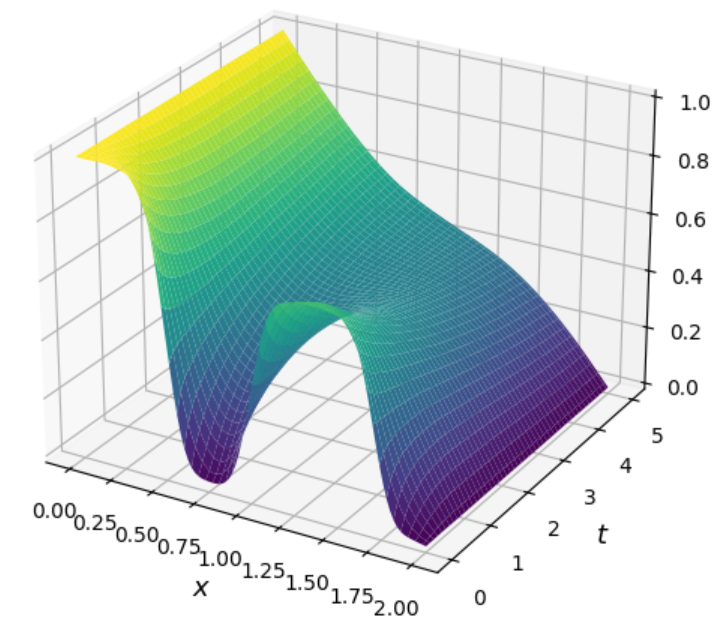}
        \caption{PINNs solution for $n_1=10$, $n_2=1$.}
        \label{fig:big n1}
    \end{minipage}\hfill
    \begin{minipage}{0.32\textwidth}
        \centering
        \includegraphics[width=\linewidth]{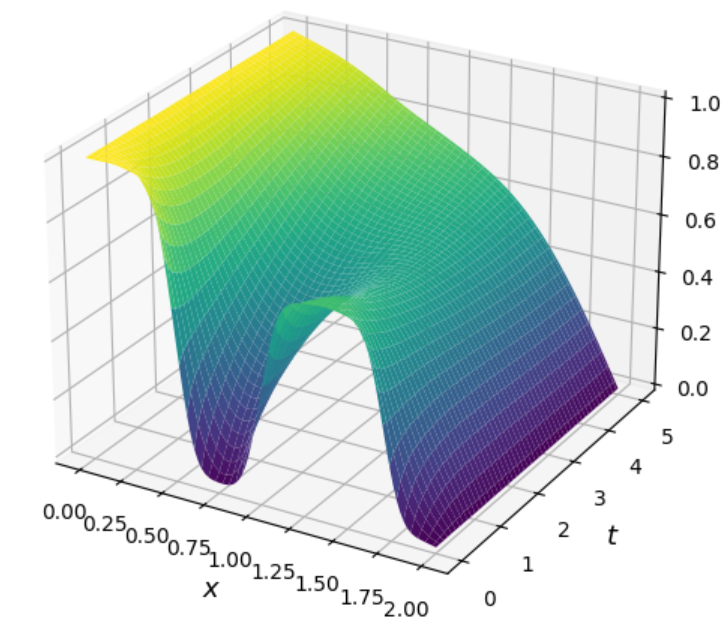}
        \caption{PINNs solution for $n_1=1$, $n_2=10$.}
        \label{fig:big n2}
    \end{minipage}\hfill
    \begin{minipage}{0.32\textwidth}
        \centering
        \includegraphics[width=\linewidth]{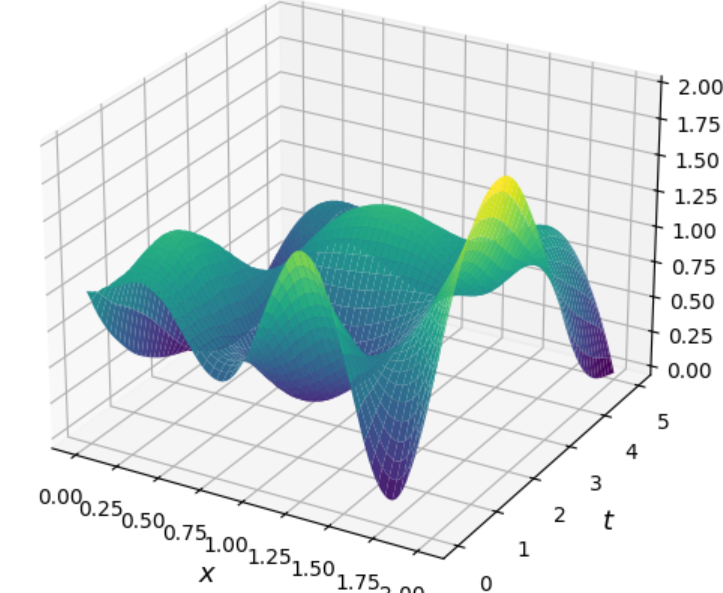}
        \caption{PINNs solution for $n_1=2$, $n_2=1$ with various conditions.}
        \label{fig:various}
    \end{minipage}
\end{figure}

\subsection{Graphs of error}
\begin{figure}[ht]
    \centering
    \begin{minipage}{0.4\textwidth}
        \centering
        \includegraphics[width=\linewidth]{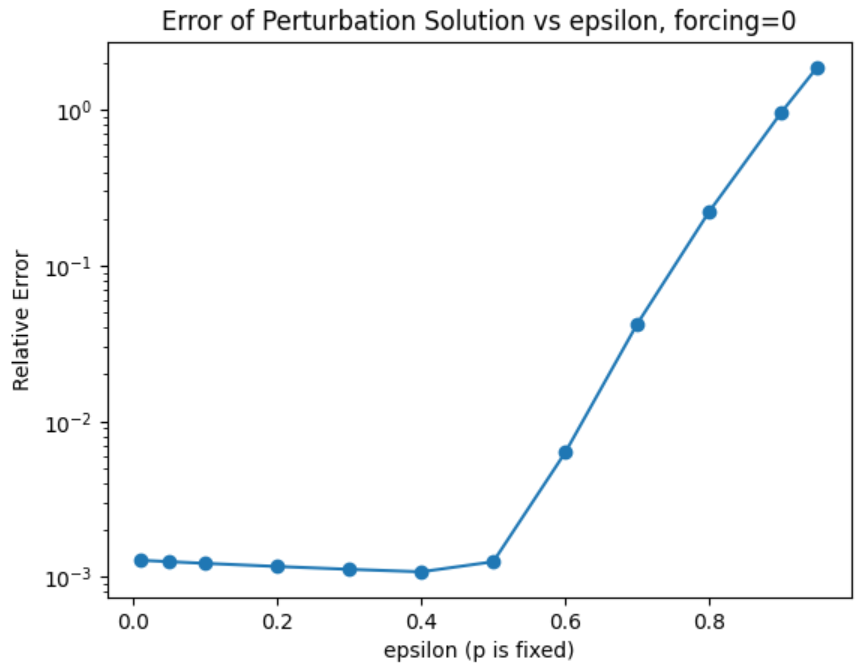}
        \caption{Relative error for different $\epsilon$, $f=0$.}
        \label{fig:epsilon error f0}
    \end{minipage}\hspace{1cm}
    \begin{minipage}{0.4\textwidth}
        \centering
        \includegraphics[width=\linewidth]{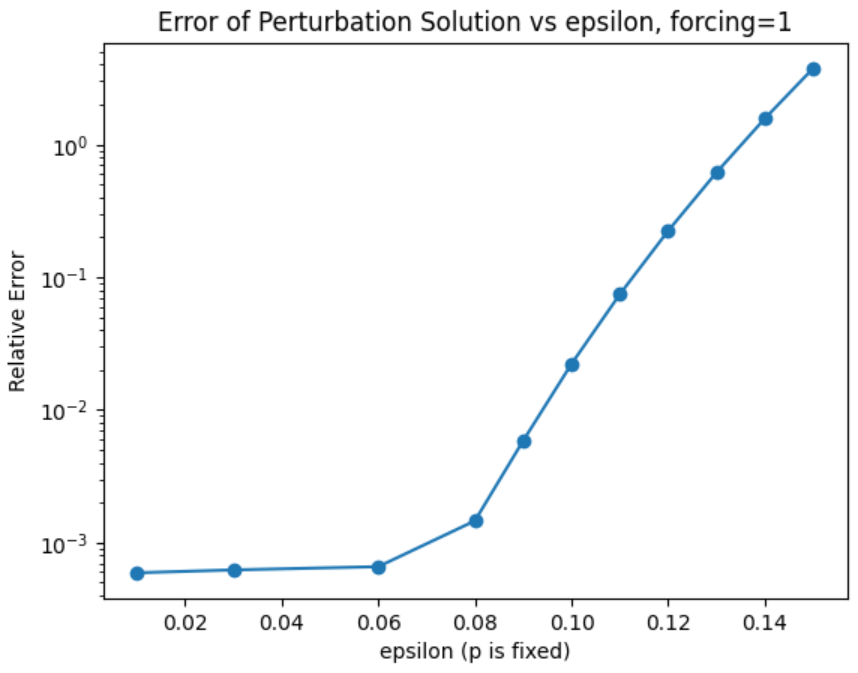}
        \caption{Relative error for different $\epsilon$, $f=1$.}
        \label{fig:epsilon error f1}
    \end{minipage}
\end{figure}

\subsection{Training Details}
\label{subsec: training}
The following diagram illustrates the overall architecture and training process:

\begin{figure}[ht]
    \centering
    \includegraphics[width=0.5\linewidth]{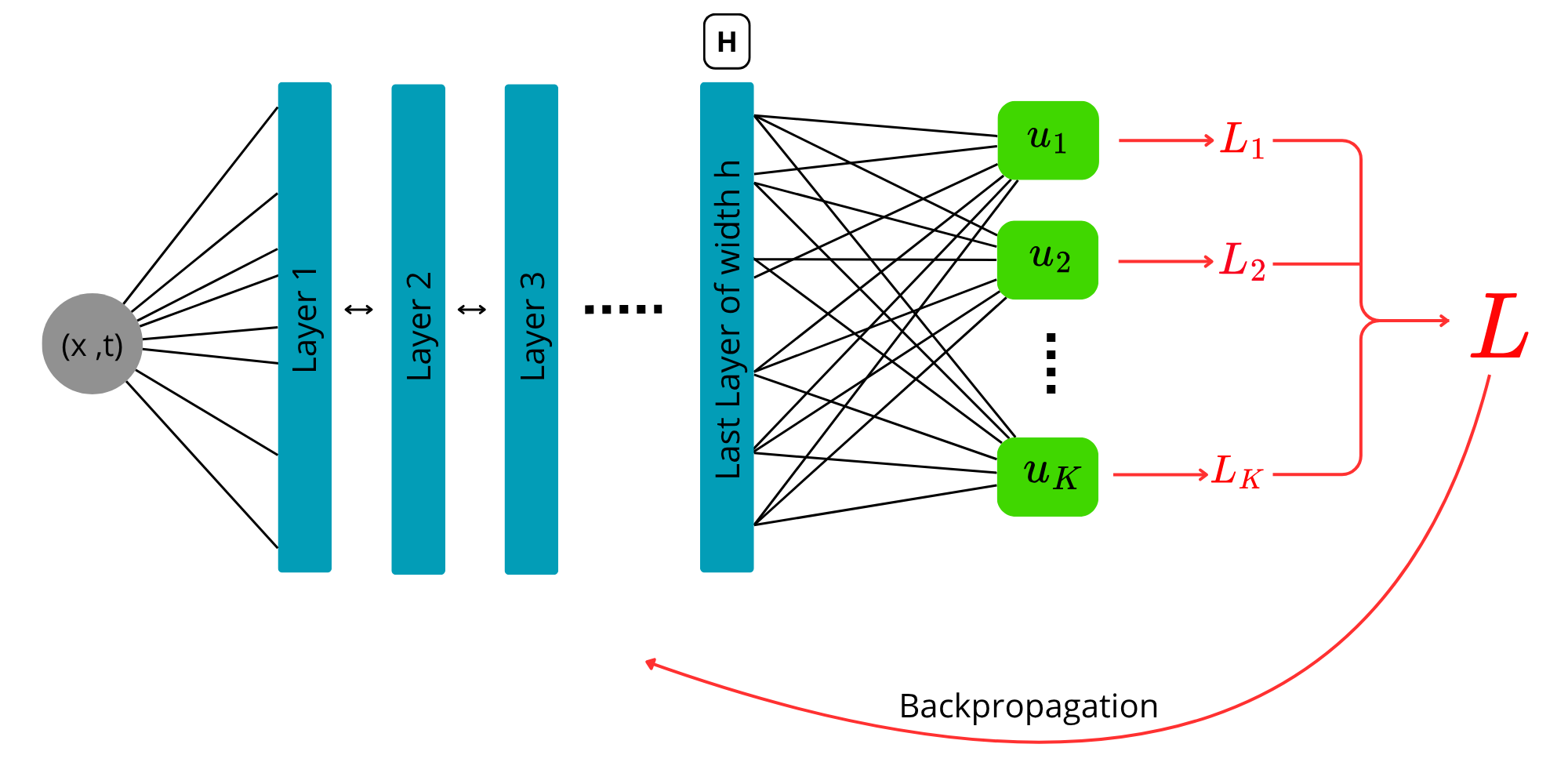}
    \caption{Architecture and process of the MH-PINN.}
    \label{fig:enter-label}
\end{figure}

The model was trained on an \textbf{Intel(R) Core(TM) i5-1035G1 CPU @ 1.00GHz} processor with \textbf{8 GB of RAM}. The network was optimized using PyTorch's Adam optimizer for 50{,}000 epochs with an initial learning rate of \( 1 \times 10^{-4} \). A scheduler reduced the learning rate by a factor of 0.975 every 1{,}000 iterations. Each iteration used 100 randomly sampled points from the domain, generated via an equally spaced grid with added noise.

\begin{figure}[ht]
    \centering
    \includegraphics[width=0.99\linewidth]{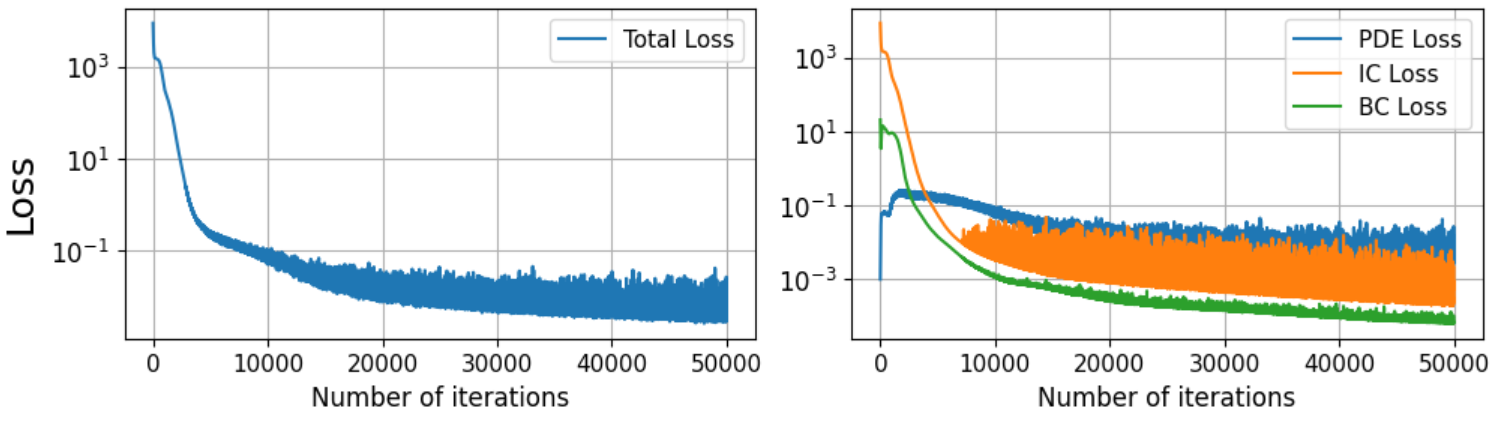}
    \caption{Total loss and components over iterations. The final total loss was \( 5.7 \times 10^{-3} \).}
    \label{fig:loss}
\end{figure}

\subsection{Results on the wave equation}
\label{Results wave equation}
For a nonlinear equation of the form
\[
\frac{\partial^2 u}{\partial t^2} - c^2\frac{\partial^2 u}{\partial x^2} + \epsilon P(u) = f(x,t)
\]
with parameters:
\begin{align*}
\text{(IC)}\ & u(x,0) = 0, \quad \frac{\partial u}{\partial t}(x,0) = 1, \\
\text{(BC)}\ & u(0,t) = 0, \quad u(2,t) = 0, \\
\text{(Forcing)}\ & f(x,t) = 0, \\
\text{(Polynomial)}\ & P(u) = u - \frac{1}{6}u^3, \quad \epsilon = 0.75, \quad c = 1,
\end{align*}
we obtain the following error versus perturbation order $p$:

\begin{figure}[ht]
    \centering
    \includegraphics[width=0.6\linewidth]{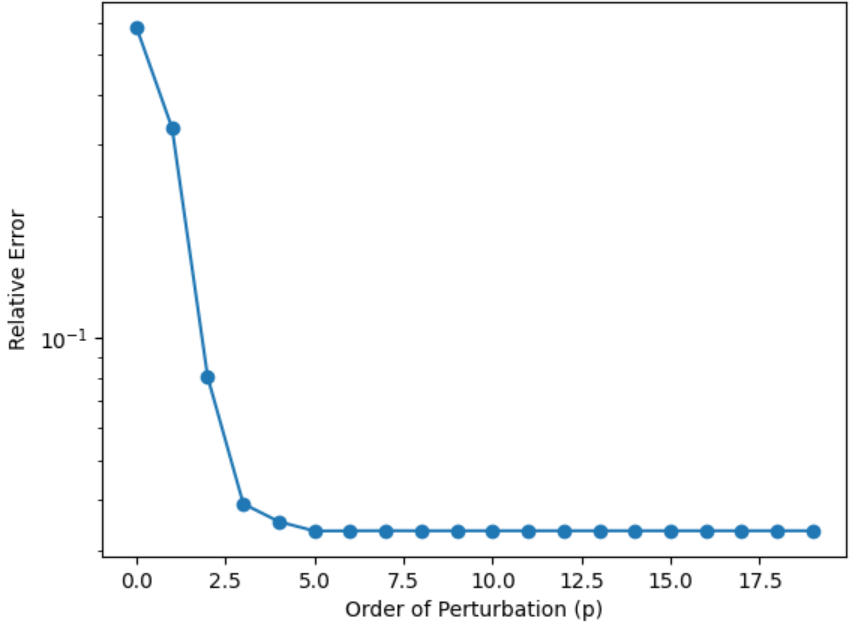}
    \caption{Relative error versus perturbation order for the modified wave equation.}
    \label{fig:perturbation wave}
\end{figure}

With the corresponding solution for $p=12$:
\begin{figure}[ht]
    \centering
    \includegraphics[width=0.75\linewidth]{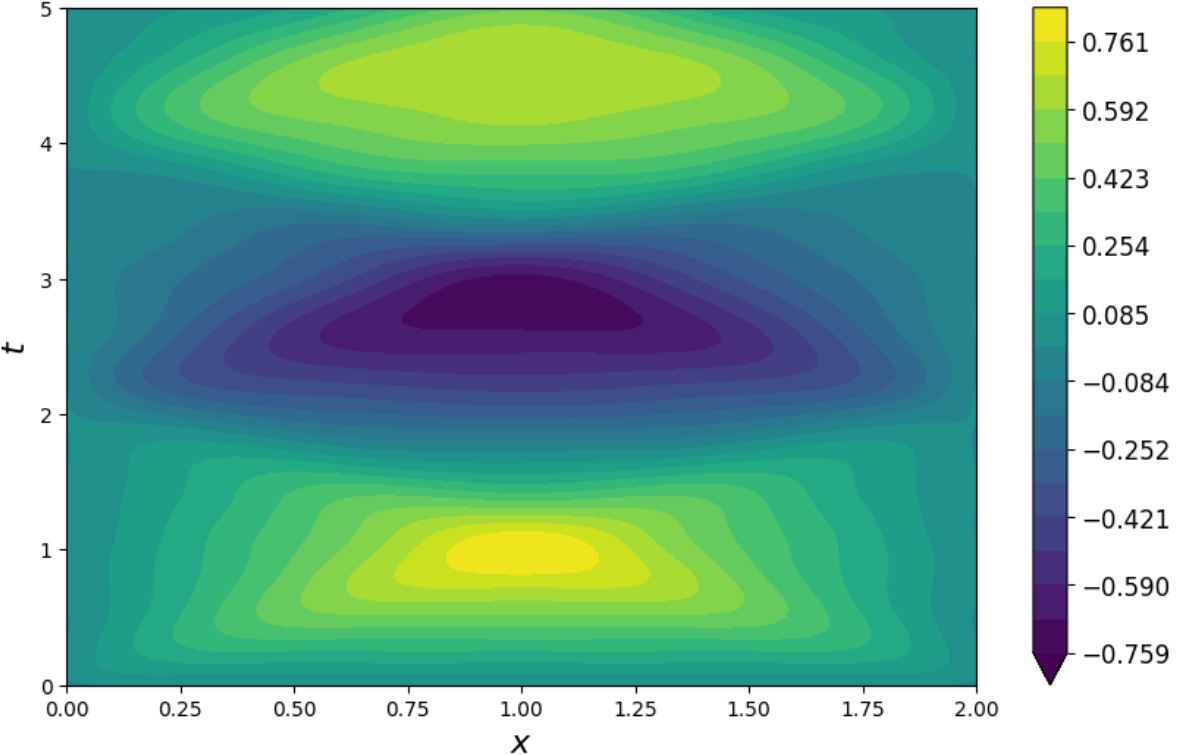}
    \caption{Solution for the perturbed wave equation.}
    \label{fig:placeholder}
\end{figure}
\newpage

\bibliographystyle{plainnat} 
\bibliography{references}

\end{document}